\documentclass{article}
\usepackage{amsthm,amsfonts,amsbsy, amssymb,amsmath,graphicx}
\usepackage{graphics}

\author{Vassily O.Manturov}
\title{Virtual links are algorithmically recognisable}
\newtheorem{thm}{Theorem}
\newtheorem{lm}{Lemma}
\newtheorem{re}{Remark}
\newtheorem{prop}{Proposition}

\begin{document} \maketitle

Virtual links were proposed by Kauffman in 1996, see \cite{Kau},
as a generalisation of classical links from the point of view of
planar diagrams.

In the present paper we prove the following

\begin{thm}
There is an algorithm to decide whether two virtual links are
equivalent or not.
\end{thm}

First of all, we use the result by Moise \cite{Moi} that all
$3$-manifoldss admit a triangulation. In the sequel, we suppose
all manifold to be triangulated.

The proof will use methods of three-dimensional topology coming
from Haken's theory of normal surfaces \cite{Ha}, and its
development by Matveev \cite{Mat}. We shall deal with
$3$-manifolds and $2$-surfaces in them. A compact surface $F$ in a
$3$-manifold  $M$ is called proper if $F\cap \partial M=\partial
F$. A proper submanifold is called {\em essential} if it does not
cobound a ball together with a part of the boundary of the
manifold. We are interested in essential spheres and essential
annuli.

A $3$-manifold $M$ is called {\em irreducible} if every $2$-sphere
in $M$ bounds a $3$-ball in $M$.

We shall use the following definition of virtual link.

A {\em virtual link} $L$ is an equivalence class for embeddings of
split union of circles into ${\cal M}\times I$ modulo isotopy,
fibre-preserving homeomorphisms of ${\cal M}\times I$ and
stabilisations/destabilisations of ${\cal M}$; we shall use the
same word ``link'' also for the image of the split union of
circles for a given realisation of the virtual link.

Here by ${\cal M}$ we mean a compact two-dimensional surface; this
surface is not necessarily connected, but for each component
${\cal M}_{i}$, the manifold ${\cal M}_{i}\times I$ should contain
at least one component of the link $L$. By
stablisations/destabilisations we mean addition/removal of
$1$-handles to ${\cal M}$ (respectively, thick cylinders to ${\cal
M}\times {I}$) such that the handle to be added/removed is
disjoint from the link.

Having a virtual link, one can try to minimise its representative
with respect to the genus of the surface ${\cal M}$ by
destablising the diagram while possible. Clearly, this can be done
in different ways. Here we should take care that the
destabilisation should be performed in such a way that no
``empty'' component appear; each $M_{i}\times I$ should contain at
least one link component. If there is a destabilisation which
divides the manifold into two components, one of which contains no
components of the link, then we should remove this component and
paste the remaining one: we attach $D^{2}\times I$ to $S^{1}\times
I$.

It is worth mentioning that {\em classical links} which are
essentially contained in a ball, admit a presentation in
$S^{2}\times I$ (or, in several copies of $S^{2}\times I$, if we
deal with a split classical link).

One of the key points that we shall use is the following theorem
due to Kuperberg, \cite{Kup}.

\begin{thm}
The result of destablisiation is well defined up to
isotopy.\label{Kpbg}
\end{thm}

Thus, in order to compare virtual knots, it is sufficient to find
their minimal representatives to and compare them.

We shall use the several theorems from Haken's theory of normal
surfaces, see \cite{Mat}.

By a {\em pattern} we mean a fixed $1$-polyhedron(graph) on the
boundary of a $3$-manifold without isolated vertices (we may
assume that this graph is a subpolyhedron of the $1$-frame of the
triangulation for the boundary). In the sequel, all proper
surfaces with boundary in a manifold $M$ with a boundary pattern
$\Gamma\subset
\partial M$ are thought to be in general position to pattern $\Gamma$, unless
otherwise specified.

A manifold $(M,\Gamma)$ with a boundary pattern is {\em boundary
irreducible} if for every proper disc $D$ the curve $\partial D$
is in general position to $\Gamma$ and bounds a disc in $\partial
M$.

By a {\em compressing disc} for a proper surface $F$ in a
$3$-manifold $M$ we mean a disc $D\subset M$, such that $D\cap
F=\partial D$. A surface $F\subset M$ is called {\em
incompressible} if for any compressing disc $D$ the curve
$\partial D$ is trivial in $D$.

A $3$-manifold is $M$ {\em sufficiently large} if there is an
incompressible closed surface $F\subset M$ which is two-sided and
different from $S^{2},{\bf RP}^{2}$.

A $3$-manifold without boundary is called {\em Haken} if it is
irreducible, boundary irreducible and sufficiently large. An
irreducible boundary irreducible $3$-manifold $(M,\Gamma)$ with
boundary-pattern $\Gamma$ is {\em Haken} if it is either
sufficiently large or $\Gamma$ is non-empty (hence, $\partial M$
is non-empty) and $M$ is a handlebody but not a ball.

Note, that a solid torus with a non-empty boundary pattern is thus
Haken.

A well-known statement \cite{Mat} says the following:

\begin{prop}
An irreducible boundary irreducible $3$-manifold with nonempty
boundary is either sufficiently large or a handlebody.\label{prp}
\end{prop}

\begin{lm}
There is an algorithm to decide whether a $3$-manifold $M$ is
reducible; if it is, the algorithm constructs a $3$-sphere
$S\subset M$ which does not bound any ball in $M$.
\end{lm}

\begin{lm}
Classical links are algorithmically recognisable.
\end{lm}

\begin{re}
The proof for an algorithm for recognising classical links, goes
in the same vein as the one we are perforiming for virtual links,
see \cite{Mat}.
\end{re}

\begin{lm}
There is an algorithm to decide whether a Haken manifold $M$ with
a pattern $\Gamma$ on the boundary $\partial M$ has an essential
proper annulus.
\end{lm}

\begin{lm}
There is an algorithm to recognise whether two Haken manifolds
$(M,\Gamma)$ and $(M',\Gamma')$ with patterns on the boundary are
homeomorphic by a homeomorphism taking $\Gamma$ to $\Gamma'$.
\end{lm}

Also, we need one more lemma.

Let $(M,L)$ be a representative of a virtual link $L$, i.e.
$M={\tilde M}\times I$ for a closed $2$-surface ${\tilde M}$. Let
$N$ be a small open tubular neighbourhood of $L$. Let us cut $N$
from $M$. We obtain a manifold to be denoted by $M_{L}$. Its
boundary consists of boundary components of $M$ (two if $M$ is
connected) and several tori; the number of tori is equal to the
number of components of $L$. Let us endow each such torus with a
pattern $\Gamma_{L}$ that represents the meridian of the
corresponding component. We obtain a manifold
$(M_{L},\Gamma_{L})$.

Clearly, virtual link $L$ can be restored from
$(M_{L},\Gamma_{L})$: we know how to make the manifold $M$ by
pasting solid tori to the boundary components of $M_{L}$ since we
know the meridians of these tori. Thus, we can restore the pair
$(M,L)$.

Suppose the link is not a split sum of a classical link with
another link.

\begin{lm}
The manifold $(M_{L},\Gamma_{L})$ is Haken.
\end{lm}

\begin{proof}
In virtue of Proposition \ref{prp}, it is sufficient to prove that
our manifold is irreducible and boundary irreducible: it never
happens for such a manifold to be a handlebody.

Now, for any compact oriented $2$-surface $S_{g}$, the manifold
$S_{g}\times I$ is irreducible, unless $g=0$. So, if our link $L$
is not classical, then for its neighbourhood $N(L)$, the set
$(S_{g}\times I)\backslash N(L)$ may be irreducible if it contains
a sphere $S$, such that $S$ bounds a ball in $S_{g}\times \{0,1\}$
containing some components of $N$. This means that those component
form a classical link which is split from the remaining components
of $L$, which leads to a contradiction.

Further, since the link $L$ is not a split sum of an unknot with a
virtual link, the manifold $M_{L}$ is boundary irreducible.

Thus, our manifold is irreducible, boundary irreducible, and hence
(by Proposition \ref{prp}) Haken.
\end{proof}

Now, let us prove the main theorem.

\begin{lm}
Let $L$ be a link that can not be reresented as a split union
$K\cup L'$ of a classical link $K$ with a virtual link $L'$. Then
$M_{L}$ is irreducible.
\end{lm}

Let $L,L'$ be a virtual links.

\begin{enumerate}

\item[Step 1.] Construct some representatives of links $(M,L),
(M',L')$. Take the corresponding manifolds with patterns on the
boundary by $(M_{L},\Gamma),(M'_{L'},\Gamma')$.

\item[Step 2.] Find whether one of $(M_{L},\Gamma)$ or
$(M'_{L'},\Gamma')$ is reducible. If yes, then it is possible to
find a reducing sphere, and, thus, extract all classical
splittable sublinks from $M_{L}$. Classical links are
algorithmically recognisable. We can compare classical split
components for $(M_{L},\Gamma)$ and $(M_{L'},\Gamma')$. If they do
not represent isotopic classical links, we stop: the initial links
are not equivalent. Otherwise, we go on.

Then, we reduce our problem to the case when no classical split
sublinks are possible. From now on, all $3$-manifolds are Haken.

\item[Step 3.]  Now, each connected component of $(M_{L},\Gamma)$
and $(M'_{L'},\Gamma')$ is a Haken manifold with a pattern on the
boundary. Thus, we can algorithmically decide whether there is a
homeomorphism $f:M_{L}\to M'_{L'}$ taking $\Gamma$ to $\Gamma'$.
If there is any, then our virtual links $L,L'$ are equivalent; if
not, they are not then $L$ and $L'$ are not equivalent virtual
links.

\end{enumerate}

Performing the  steps described above, we obtain the claim of the
theorem.

\begin{re}
The proof described above works equivalently for oriented virtual
links and for framed virtual links.
\end{re}

\end{document}